\documentclass[12pt,a4paper]{article}

\usepackage{t1enc}
\usepackage[leqno]{amsmath}
\usepackage{amssymb,latexsym}
\usepackage{pictex}
\usepackage{theorem}

\newtheorem{theorem}{Theorem}[section]
\newtheorem{proposition}[theorem]{Proposition}
\newtheorem{lemma}[theorem]{Lemma}
\newtheorem{corollary}[theorem]{Corollary}
\newtheorem{definition}[theorem]{Definition}

\newenvironment{proof}{\par\addvspace{12pt plus 5pt minus 3pt}\noindent\textsc{Proof: }}{\qed\par\addvspace{8pt plus 3pt minus 1.5pt}}

\newcommand{\ov}{\overline}
\newcommand{\K}{\mathcal{K}}
\newcommand{\Knice}{\K_{\mathrm{nice}}}

\newcommand{\id}{{\mathrm{id}}}
\renewcommand{\L}{\mathcal{L}}
\renewcommand{\P}{\mathcal{P}}
\renewcommand{\S}{\mathcal{S}}
\newcommand{\T}{\mathcal{T}}
\newcommand{\U}{\mathcal{U}}
\newcommand{\V}{\mathcal{V}}
\newcommand{\W}{\mathcal{W}}
\newcommand{\X}{\mathcal{X}}
\newcommand{\chara}{\mathrm{char}}
\renewcommand{\cot}{\mathrm{cot}}
\newcommand{\coi}{\mathrm{coi}}
\newcommand{\cof}{\mathrm{cof}}
\newcommand{\aut}{\mathrm{Aut}}
\newcommand{\out}{\mathrm{Out}}
\newcommand{\inn}{\mathrm{Inn}}
\newcommand{\sym}{\mathrm{Sym}}

\newcommand{\qed}{\hfill $\Box$}

\begin{document}
\title{Outer Automorphism Groups of Ordered Permutation Groups}
\author{Manfred Droste\\
  \small Institut f\"ur Algebra,
  Technische Universit\"at Dresden,\\[-1ex]
  \small D-01062 Dresden, Germany\\[-1ex]
  \small e-mail: droste@math.tu-dresden.de
  \and
  Saharon Shelah\thanks{Research supported by the German-Israeli-Foundation for
  Scientific Research and Development. Publication 743.}\\
  \small Institute of Mathematics, Hebrew University,\\[-1ex]
  \small Jerusalem, Israel\\[-1ex]
  \small e-mail: shelah@math.huji.ac.il
}
\date{}
\maketitle

\begin{abstract}
An infinite linearly ordered set $(S,\le)$ is called doubly homogeneous if its
automorphism group $A(S)$ acts $2$-transitively on it. We show that any group
$G$ arises as outer automorphism group $G\cong\out(A(S))$ of the automorphism
group $A(S)$, for some doubly homogeneous chain $(S,\le)$.
\end{abstract}

\section{Introduction}

An infinite linearly ordered set (``chain'') $(S,\le)$ is called \emph{doubly
  homogeneous}, if its automorphism group, i.e.~the group of all
order-preserving permutations, $A(S)=\aut(S,\le)$ acts $2$-transitively on
it. Chains $(S,\le)$ of this type and their automorphism groups $A(S)$ have
been intensively studied. They have been used e.g.~for the construction of
infinite simple torsion free groups (Higman~\cite{Hi}) or, in the theory of
lattice-ordered groups ($\ell$-groups), for embedding arbitrary $\ell$-groups
into simple divisible $\ell$-groups (Holland~\cite{H1}). The normal subgroup
lattices of the groups $A(S)$ have been determined in \cite{BD,DS}. Obviously,
all linearly ordered fields are examples for such chains. For a variety of
further results, see Glass~\cite{G}. Here, we will be concerned with outer
automorphism groups $\out(A(S))=\aut(A(S))/\inn(A(S))$ of the automorphism
groups $A(S)$ for doubly homogeneous chains $(S,\le)$.

In the literature, many authors have dealt with the problem of determining
which groups $G$ can arise as $G\cong\out(H)$ for $H$ in a given class of
groups. Already Schreier and Ulam~\cite{SU} showed that each automorphism of
the infinite symmetric groups is inner,
i.e.~$\out(\sym(\mathbb{N}))=\{1\}$. It is known that each outer automorphism
of $A(\mathbb{Q})$ or of $A(\mathbb{R})$ arises via conjugation by an
anti-automorphism of $(\mathbb{R},\le)$, hence
$\out(A(\mathbb{Q}))\cong\out(A(\mathbb{R}))\cong\mathbb{Z}_2$, and it is
known how to construct doubly homogeneous chains $(S,\le)$ for which
$\out(A(S))$ is trivial (Holland~\cite{H2}, Weinberg~\cite{W},
Droste~\cite{D}). Solving a problem which had been open for quite some time,
Holland~\cite{H2} constructed a doubly homogeneous chain 
$(S,\le)$ for which $A(S)$ has an outer automorphism \emph{not} arising from an
anti-automorphism of the Dedekind-completion $(\overline{S},\le)$ of
$(S,\le)$. In his example, $\out(A(S))\cong V_4$, Klein's four-group. Assuming
the generalized continuum hypothesis (GCH), McCleary~\cite{Mc2} constructed
further examples of this type. However, to date $V_4$, $\mathbb{Z}_2$
and the trivial group are the only groups realized as $\out(A(S))$, where
$(S,\le)$ is a doubly homogeneous chain. In fact, a realization of any group
$G$ even just as the outer automorphism group $G\cong\out(H)$ of some group $H$
was established more recently in Matumoto~\cite{Ma}.

With the pointwise ordering of functions, $A(S)$ becomes an $\ell$-group, and
Holland~\cite{H1} showed that any $\ell$-group $H$ can be $\ell$-embedded
(i.e.~embedded as an $\ell$-group) into $A(S)$, for some doubly homogeneous
chain $(S,\le)$. Here, we will show the following generalization of the
previously mentioned results:

\begin{theorem}\label{theorem main}
Let $G$ be any group, $H$ any $\ell$-group and $\lambda$ a regular uncountable
cardinal with $\lambda\ge\lvert G\rvert$ and $\lambda>\lvert H\rvert$. Then
there exists a doubly homogeneous chain $(S,\le)$ of cardinality $\lambda$
such that $G\cong\out(A(S))$ and $H$ $\ell$-embeds into $A(S)$.
\end{theorem}

Here, the realization result $G\cong\out(A(S))$ involves constructing a doubly
homogeneous chain $(S,\le)$ such that $A(\overline{S})$ acts on the set of
orbits of $A(S)$ in $\overline{S}$ just like $G$. Using codings of the group
action of $G$ through a system of suitable stationary subsets of $\lambda$
inside $\overline{S}$, we will first describe a class of doubly homogeneous
chains $(S,\le)$ for which $G\cong\out(A(S))$ follows. Then, in
Section~\ref{sec 4}, we will actually construct these chains. As often for
homogeneous structures, this could be done by suitable amalgamations of linear
orderings, but here we will use methods of \cite{DS} for a more explicit
construction. For the simultaneous embedding of $H$ into $A(S)$, we will use
Holland's result \cite{H1}.

By a L\"owenheim-Skolem argument, we obtain as a consequence:

\begin{corollary}\label{cor main}
Let $H$ be any $\ell$-group, $G$ any group and $\lambda$ a regular uncountable
cardinal with $\lambda\ge\lvert G\rvert$ and $\lambda>\lvert H\rvert$. Then
there exists an $\ell$-group $K$ with $H\subseteq K$ (as $\ell$-groups) and $\lvert
K\rvert=\lambda$ such that $\out(K)\cong G$.
\end{corollary}

Here, each group automorphism of $K$ is also a lattice automorphism of (the
$\ell$-group) $K$.

Similar realization results have been established in the literature for
various classes of groups. Dugas and G\"obel~\cite{DG1} showed that any
countable group
arises as the outer automorphism group of some locally-finite $p$-group, and
G\"obel and Paras~\cite{GP} established the corresponding result for the class
of torsion-free metabelian groups. Note that the groups $A(S)$ of
Theorem~\ref{theorem main} are not simple. In Droste, Giraudet and
G\"obel~\cite{DGG}, Theorem~\ref{theorem main} will be used to show that \emph{any
group can be realized as outer automorphism group of a simple group}, which in
turn arises as the automorphism group of a suitable homogeneous circle.

\section{Outer automorphisms}
\label{sec 2}

In this section, we will present a class of doubly homogeneous chains
$(S,\le)$ for which $\out(A(S))\cong G$, for a given group $G$.

For any chain $(S,\le)$, we denote by $\overline{S}$ its
Dedekind-completion. Clearly, each automorphism $f$ of $(S,\le)$ extends
uniquely to an automorphism of $(\overline{S},\le)$ which we will also denote
by $f$; hence $A(S)\subseteq A(\overline{S})$. We recall that $(S,\le)$ is
\emph{doubly homogeneous}, if for all $u,v,x,y\in S$ with $u<v$ and $x<y$
there is $g\in A(S)$ such that $u^g=x$ and $v^g=y$.
If $x\in\overline{S}$, the set $x^{A(S)}=\{x^f:\ f\in A(S)\}$ is called an
\emph{orbit} of $A(S)$ in $\overline{S}$. The following result relates outer
automorphisms of $A(S)$ to automorphisms of $\overline{S}$ permuting the
orbits of $A(S)$.

\begin{proposition}[\protect{\cite{G}}]\label{prop 2.1}
Let $(S,\le)$ be a doubly homogeneous chain. Each automorphism $\varphi$ of
$A(S)$ corresponds bijectively to an automorphism or anti-automorphism $f$ of
$(\overline{S},\le)$, which permutes the orbits of $A(S)$ in $\overline{S}$,
such that
$g^\varphi=f^{-1}\circ g\circ f\ \text{ for all }g\in A(S)$.
\end{proposition}

Now assume $(\overline{S},\le)$ is not anti-isomorphic to itself, and let
$Z:=\{U\subseteq\overline{S}:\ U \text{ orbit of } A(S) \text{ in }
\overline{S}, \text{ and } (S,\le)\cong (U,\le)\}$. By Proposition~\ref{prop
  2.1}, each automorphism $\varphi$ of $A(S)$ induces an element $f\in
A(\overline{S})$ which permutes $Z$, and we obtain an epimorphism from
$\aut(A(S))$ onto the group $G_Z=\{p\in \sym(Z):\ \exists f\in
A(\overline{S}).~f \text{ induces } p\}$ with kernel $\inn(A(S))$. Hence
$\out(A(S))\cong G_Z$. In order to prove Theorem~\ref{theorem main}, given any
group $G$ we therefore have to construct a doubly homogeneous chain $(S,\le)$
with $(\overline{S},\le)$ not anti-isomorphic to itself such that $G$
represents the action of $A(\overline{S})$ on the orbits of $A(S)$ in
$\overline{S}$, i.e.~$G\cong G_Z$.

We recall some notation, which is mostly standard. Let $(S,\le)$ be
any chain. If $x\in\overline{S}$ has no immediate predecessor, we define the
\emph{cofinality} of $x$ to be
\[
\cof(x):=\min\{\lvert A\rvert:\ A\subseteq\overline{S}, x\not\in A, x=\sup
A\}.
\]
We adopt the convention that if $x$ has an immediate predecessor, then
$\cof(x)=\aleph_0$. We define the \emph{coinitiality} $\coi(x)$ of $x$
dually. If $\cof(x)=\coi(x)$ then this cardinal is called the
\emph{coterminality} of $x$, denoted $\cot(x)$. The ordered pair
$(\cof(x),\coi(x))$ is called the \emph{character} of $x$, denoted $\chara(x)$,
and we let $\chara(S)=\{\chara(x):\ x\in S\}$. If $X\subseteq\overline{S}$, we
let $\chara_{\overline{S}}(X)=\{\chara(x):\ x\in X\}$, where $\chara(x)$ is
determined in $\overline{S}$. If $a,b\in\overline{S}$ with $a<b$
and $X\subseteq\overline{S}$, let $[a,b]_X=\{x\in X:\ a\le x\le b\}$. If
$A,B\subseteq\overline{S}$, we write $A<B$ to denote that $a<b$ for all $a\in
A$, $b\in B$, and $A<x$ abbreviates $A<\{x\}$. The
chain $(S,\le)$ has \emph{countable coterminality}, denoted $\cot(S)=\aleph_0$,
if it contains a countable subset which is unbounded above and below in
$S$. We say that $(S,\le)$ is \emph{dense}, if for all $a<b$ in $S$ there is
$s\in S$ with $a<s<b$, and \emph{unbounded}, if $S$ does not contain a
greatest or smallest element.

Clearly, all points in a given orbit of $A(S)$ in $\overline{S}$ have the same
character. Also, it is well-known and easy to see by piecewise patching of
automorphisms of $S$ (cf.~\cite{G} or below the argument for Lemma~\ref{lemma
  3.1}), that if $(S,\le)$ is doubly homogeneous, then all elements
$x\in\overline{S}\setminus S$ with $\cot(x)=\aleph_0$ form a single orbit of
$A(S)$ in $\overline{S}$. Hence when constructing the chain $(S,\le)$ for
Theorem~\ref{theorem main}, we have to resort to orbits of $A(S)$ whose
elements do not have countable coterminality.

As usual, we identify cardinals with the least ordinal of their cardinality. A
subset $A$ of a cardinal $\lambda$ is called \emph{stationary}, if $A\cap
C\not=\emptyset$ for each closed unbounded subset $C$ of $\lambda$. A sequence
$\langle x_i:\ i\in\lambda\rangle\subseteq\overline{S}$ is called
\emph{continuously increasing} if $x_i<x_j$ for each $i<j<\lambda$, and
$x_j=\sup_{i<j}x_i$ in $\overline{S}$ for each limit ordinal $j<\lambda$. We
define \emph{continuously decreasing} dually.

Now we make the following

\paragraph{General assumption.}
In all of this section let $G$ be a group with neutral element $e$, and let
$\aleph_0\le\kappa<\lambda$ be two regular cardinals such that
$\lambda\ge\lvert G\rvert$. Let $S^\lambda_\kappa=\{j<\lambda:\
\cof(j)=\kappa\}$, a stationary subset of $\lambda$. By Solovay's theorem
\cite{S}, we split
$S^\lambda_\kappa=\dot{\raisebox{0pt}[0pt]{$\bigcup$}}_{g\in G}S_g$ into
$\lvert G\rvert$ pairwise disjoint stationary subsets $S_g$ ($g\in G$).

\paragraph{}
Now we define a class of chains which will be crucial in all of our subsequent
considerations.

\begin{definition}\label{def 2.2}
Let $\K$ be the class of all structures $\S=(S,\le,(P_g)_{g\in G},(Q_g)_{g\in
  G})$ with the following properties:
\begin{enumerate}
\item $(S,\le)$ is a dense unbounded chain;
\item $P_g$, $Q_g$ ($g\in G$) are pairwise disjoint non-empty subsets of
  $\overline{S}\setminus S$; let $P=\dot{\raisebox{0pt}[0pt]{$\bigcup$}}_{g\in
  G}P_g$ and $Q=\dot{\raisebox{0pt}[0pt]{$\bigcup$}}_{g\in G}Q_g$;
\item whenever $x\in P_g$ ($g\in G$), then $\chara(x)=(\lambda,\kappa)$, and
  there is a continuously increasing sequence $\langle x_i:\ i<\lambda\rangle$
  in $\overline{S}$ such that $x=\sup_{i<\lambda}x_i$ and whenever $h\in G$
  and $j\in S_h$, then $x_j\in Q_{hg}$;
\item whenever $y\in Q_g$ ($g\in G$), then $\chara(y)=(\kappa,\lambda)$, and
  there is a continuously decreasing sequence $\langle y_i:\ i<\lambda\rangle$
  in $\overline{S}$ such that $y=\inf_{i<\lambda}y_i$ and whenever $h\in G$
  and $j\in S_h$, then $y_j\in P_{hg}$;
\item whenever $a,b,c,d\in\overline{S}$ with $a<b$, $c<d$, $x\in P\cap(a,b)$,
  $y\in Q\cap(a,b)$ and $f:[a,b]_{\overline{S}}\longrightarrow
  [c,d]_{\overline{S}}$ is an order-isomorphism, then $x^f\in P$ and $y^f\in
  Q$.
\end{enumerate}
\end{definition}

We will refer to the sets $P_g$, $Q_g$ ($g\in G$) as the \emph{colours of
  $\S$}, and the sets $P_g$ ($g\in G$) are the \emph{$P$-colours}.
Condition~\ref{def 2.2}(5) means that we can ``locally recognize'' points
belonging to $P$ or to $Q$. It ensures that then
\begin{equation}\tag{$*$}
P^f\subseteq P \text{ and } Q^f\subseteq Q \text{ for each } f\in A(\overline{S}).
\end{equation}
If $\aleph_0<\kappa<\lambda$, we note that we will construct our structures
$\S$ such that
\begin{align} %\tag{$+$}
\chara_{\overline{S}}(P) &= \{x\in\overline{S}:\ \chara(x)=(\lambda,\kappa)\}
  \text{ and}\\
\chara_{\overline{S}}(Q) &= \{y\in\overline{S}:\ \chara(y)=(\kappa,\lambda)\};
\end{align}
this obviously implies condition~\ref{def 2.2}(5). Next we note:

\begin{lemma}\label{lemma 2.3}
Let $\S=(S,\le,(P_g)_{g\in G},(Q_g)_{g\in G})\in\K$. Let $f\in
\aut(\overline{S},\le)$ satisfy $P_e^f\subseteq P_g$ for some $g\in G$. Then
$P_h^f\subseteq P_{hg}$ and $Q_h^f\subseteq Q_{hg}$ for each $h\in G$. In
particular, if $P_e^f\subseteq P_e$, then each $P_h$ and each $Q_h$ ($h\in G$)
is invariant under $f$.
\end{lemma}

\begin{proof}
Let $h\in G$. We first show that $Q_h^f\subseteq Q_{hg}$. Choose any $x\in
Q_h$. By requirement \ref{def 2.2}(5), clearly $x^f\in
Q$. So, $x^f\in Q_k$ for some $k\in G$. By condition \ref{def 2.2}(4), we can
find two continuously decreasing sequences $\langle y_i:\ i<\lambda\rangle$
and $\langle z_i:\ i<\lambda\rangle$ in $\overline{S}$ such that
$x=\inf_{i<\lambda}y_i$, $x^f=\inf_{i<\lambda}z_i$, and whenever $h'\in G$ and
$j\in S_{h'}$, then $y_j\in P_{h'h}$ and $z_j\in P_{h'k}$. Standard
back-and-forth arguments involving closed unbounded subsets of uncountable
cardinals show that $C=\{i<\lambda:\ y_i^f=z_i\}$ is a closed unbounded subset
of $\lambda$. Since $S_{h^{-1}}$ is stationary, there is $j\in S_{h^{-1}}$
with $y_j^f=z_j$. So, $y_j\in P_e$ and $z_j=y_j^f\in P_g\cap
P_{h^{-1}k}$. Thus $k=hg$ as claimed.

By a similar argument, it now follows that $P_h^f\subseteq P_{hg}$ for each
$h\in G$. The final statement is then immediate.
\end{proof}

Now let $\S=(S,\le,(P_g)_{g\in G},(Q_g)_{g\in G})$ and
$\S'=(S',\le,(P'_g)_{g\in G},(Q'_g)_{g\in G})\in\K$ and $A\subseteq
\overline{S}$, $B\subseteq\overline{S'}$. Let $h\in G$ and $f:
(A,\le)\longrightarrow(B,\le)$ be an order-isomorphism. We say that $f$
\emph{maps the colours as prescribed by} $h$, if
\[
(A\cap P_g)^f=B\cap P'_{gh} \text{ and } (A\cap Q_g)^f=B\cap Q'_{gh}\quad \text{for each $g\in G$.}
\]
If here $\S=\S'$, we also say that $f$ \emph{permutes the
  colours as prescribed by} $h$.

We call $\S$ \emph{doubly homogeneous in the $P$-colours}, if for any $g,h\in
G$ and $u,v\in P_g$, $x,y\in P_{gh}$ with $u<v$ and $x<y$, there is $f\in
A(S)$ such that $\{u,v\}^f=\{x,y\}$ and $f$ permutes the colours as prescribed
by $h$.

Now we can prove:

\begin{theorem}\label{theorem 2.4}
Let $\S=(S,\le,(P_g)_{g\in G},(Q_g)_{g\in G})\in\K$ be doubly homogeneous in
the $P$-colours such that $(\overline{S},\le)$ is not anti-isomorphic to
itself, and let $X=P_e$. Then $(X,\le)$ is a doubly homogeneous chain and
$\out(A(X))\cong G$.
\end{theorem}

\begin{proof}
Clearly, $(X,\le)$ is doubly homogeneous and dense in $\overline{S}$. We
identify $\overline{X}=\overline{S}$. By Lemma~\ref{lemma 2.3}, each $P_g$ and
each $Q_g$ ($g\in G$) is invariant under $A(X)$. By homogeneity we get
$P_g=P_g^{A(X)}$ for each $g\in G$. Hence each $P_g$ ($g\in G$) is an orbit of
$A(X)$ in $\overline{X}$. By Proposition~\ref{prop 2.1}, any outer
automorphism $\psi$ of $A(X)$ determines an automorphism $f$ of $\overline{X}$
which permutes the orbits of $A(X)$ in $\overline{X}$ and hence by
requirement \ref{def 2.2}(3), permutes the orbits $P_g$ ($g\in G$) among
themselves. By Lemma~\ref{lemma 2.3}, this permutation $f$ determines an
element $h$ of $G$. Conversely, by homogeneity of $\S$, any element $h$ of $G$
can be realized in this way by an automorphism $f$ of $\overline{X}$ permuting
the sets $P_g$ ($g\in G$), and hence $h$ is realized by an outer automorphism
$\psi$ of $A(X)$. This correspondence constitutes the required isomorphism.
\end{proof}

\section{Isomorphisms between intervals}
\label{sec 3}

It is well-known that a chain $(S,\le)$ is doubly homogeneous if and only if
any two of its intervals $[a,b]_S$ and $[c,d]_S$ ($a,b,c,d\in S$ with $a<b$,
$c<d$) are order-isomorphic. In this section we will derive a similar result
for particular structures in $\K$. In all of this section, we make the
\emph{general assumption} of Section~\ref{sec 2}.

Let $\S=(S,\le,(P_g)_{g\in G},(Q_g)_{g\in G})\in\K$. We say
that $\S$ is \emph{G-homogeneous for $S$-intervals}, if for any $u<v$ and $x<y$
in $S$ and any $h\in G$, there is an isomorphism $f:
[u,v]_S\longrightarrow[x,y]_S$ permuting the colours as prescribed by $h$.

We call $\S$ \emph{G-homogeneous for $P$-intervals}, if for any $g,h\in G$ and
$u,v\in P_g$, $x,y\in P_{gh}$ with $u<v$ and $x<y$, there is an isomorphism $f:
[u,v]_S\longrightarrow[x,y]_S$ permuting the colours as prescribed by $h$.

Clearly, if $\S$ is doubly homogeneous in the $P$-colours, $\S$ is also
$G$-homoge\-neous for $P$-intervals. For the converse we need additional
assumptions on $\S$. In general, there seems to be no relationship between
$G$-homogeneity for $S$-intervals and for $P$-intervals,
respectively. However, in this section we will show for particular
structures $\S\in\K$ that $G$-homogeneity for $S$-intervals implies
$G$-homogeneity for $P$-intervals. First we note:

\begin{lemma}\label{lemma 3.1}
Let $\S=(S,\le,(P_g)_{g\in G},(Q_g)_{g\in G})\in\K$ be G-homogeneous for
$S$-intervals. Let $u,v,x,y\in\overline{S}$ such that $u<v$, $x<y$ and
$\cot(u,v)=\aleph_0=\cot(x,y)$. Let $h\in G$. Then there is an isomorphism $f:
(u,v)_S\longrightarrow(x,y)_S$ permuting the colours as prescribed by $h$.
\end{lemma}

\begin{proof}
Choose $\mathbb{Z}$-sequences $(a_i)_{i\in\mathbb{Z}}\subseteq(u,v)_S$ and
$(b_i)_{i\in\mathbb{Z}}\subseteq(x,y)_S$ such that $a_i<a_{i+1}$ and
$b_i<b_{i+1}$ for each $i\in\mathbb{Z}$ and $u=\inf_{i\in\mathbb{Z}}a_i$,
$v=\sup_{i\in\mathbb{Z}}a_i$, $x=\inf_{i\in\mathbb{Z}}b_i$,
$y=\sup_{i\in\mathbb{Z}}b_i$. For each $i\in\mathbb{Z}$, there is an
isomorphism $f_i: [a_i,a_{i+1}]_S\longrightarrow[b_i,b_{i+1}]_S$ permuting the
colours as prescribed by $h$. Patching the $f_i$'s together, we obtain the
required isomorphism $f$.
\end{proof}

Now we turn to the definition of structures $\S\in\K$ for which there is a
closer relationship between the three types of homogeneity we defined.

A partially ordered set $(T,\preceq)$ is called a \emph{tree}, if it contains
a smallest element and for each $x\in T$, the set $(\{t\in T:\ t\prec
x\},\preceq)$ is well-ordered; the cardinality of this set will also be called
the \emph{height} of $x$. We write $M_x$ for the set of immediate successors
( = minimal strict upper bounds) of $x$ in $(T,\preceq)$. We will now consider
particular trees together with a partial ordering which will later on be
extended to a linear ordering.

\begin{definition}
Let $\T$ be the class of all structures $(T,\preceq,\le)$ with the following
properties:
\begin{enumerate}
\item $(T,\preceq)$ is a non-singleton tree and $\le$ is a partial order on
  $T$;
\item each element $x\in T$ has finite height and if the height of $x$ is even
  (odd), then $(M_x,\le)$ is either empty or a chain isomorphic
  (anti-isomorphic, respectively) to $\lambda$;
\item whenever $y\in M_x$ for some $x\in T$ of even (odd) height, then
  $M_y\not=\emptyset$ if and only if in $(M_x,\le)$ we have $\cof(y)=\kappa$
  ($\coi(y)=\kappa$, respectively).
\end{enumerate}
\end{definition}
Observe that each such tree has size $\lambda$.

Let $(T,\preceq,\le)\in\T$. We denote the smallest element of $(T,\preceq)$ by
$t_0$. Now we extend the partial order $\le$ to a linear order $\le'$
on $T$ as follows. First, we put $M_{t_0}<'t_0$. Now assume $x\in T$ has even
(odd, respectively) height and $y\in M_x$ with $M_y\not=\emptyset$. Then put
$y<'M_y<'y^+$ ($y^+<'M_y<'y$) where $y^+$ denotes the immediate successor
(predecessor) of $y$ in $(M_x,\le)$, respectively. For simplicity, we denote
this extension $\le'$ also by $\le$.

Observe that in the Dedekind-completion $(\overline{T},\le)$, each element
$x\in\overline{T}\setminus T$ corresponds uniquely to a maximal path
$(a_i)_{i\in\omega}$ in $(T,\le)$ such that
%%%%
\linebreak
 $a_0=t_0$, $a_{i+1}\in M_{a_i}$,
and $a_{2i+1}<a_{2i+3}<a_{2i+2}<a_{2i}$ for each $i\in\omega$, and
$x=\sup_{i\in\omega}a_{2i+1}=\inf_{i\in\omega}a_{2i}$ in
$(\overline{T},\le)$. In particular, $\cot(x)=\aleph_0$. Hence
$\chara(\overline{T})=\{(\mu,\aleph_0), (\aleph_0,\mu):\
\aleph_0\le\mu\le\lambda,\ \mu\not=\kappa, \text{ $\mu$
  regular}\}\cup
%%%%
\linebreak
\{(\lambda,\kappa), (\kappa,\lambda), (\aleph_0,\aleph_0)\}$.

Next, we define for each $k\in G$ a \emph{coloured linear ordering}
$\L_k=
%%%%
\linebreak
(L_k,\le,(P_g)_{g\in G}, (Q_g)_{g\in G})$ as follows. Let $(A,\le)$ be a
fixed linear ordering anti-isomorphic to $\kappa$ and put, irregardless of
$k$, $(L_k,\le)=(\overline{T},\le)+(A,\le)$, the disjoint sum with
$(\overline{T},\le)$ ``to the left'' of $(A,\le)$. Thus $t_0=\inf A$ in
$(L_k,\le)$. We define the sets $P_g,Q_g\subseteq T$ now. Only elements of
even (odd) height will belong to $\bigcup_{g\in G}P_g$ ($\bigcup_{g\in
  G}Q_g$), respectively. First, put $t_0\in P_k$. By induction, assume $t\in
T$ satisfies $t\in P_g$ and $(M_t,\le)\cong\lambda$, thus $M_t=\{x_i:\
i\in\lambda\}$ with $x_i<x_j$ if $i<j$ in $\lambda$. Then for each $h\in G$
and $j\in S_h$ put $x_j\in Q_{hg}$. Next, let $t\in Q_g$ and $(M_t,\le)$ be
anti-isomorphic to $\lambda$, thus $M_t=\{y_i:\ i\in\lambda\}$ with $y_j<y_i$
if $i<j$ in $\lambda$. Then for each $h\in G$ and $j\in S_h$ put $y_j\in
P_{hg}$.

The reader will notice the similarity of $\L_k$ to the structures in $\K$. In
fact, the chains $\L_k$ will be basic building blocks for particular chains in
$\K$. Observe that $|P_g|=|Q_g|=\lambda$ for each $g\in G$, since any
stationary subset of $\lambda$ has size $\lambda$. First we note:

\begin{lemma}\label{lemma 3.3}
Let $h,k\in G$.
\begin{itemize}
\item[(a)] Let $\L_h=(L_h,\le,(P_g)_{g\in G},(Q_g)_{g\in G})$ and
  $\L_{kh}=(L_{kh},\le,(P'_g)_{g\in G},(Q'_g)_{g\in G})$. Then $f=id:
  L_h\longrightarrow L_{kh}$ satisfies $P_g^f=P'_{gh}$ and $Q_g^f=Q'_{gh}$ for
  each $g\in G$.
\item[(b)] Let $\S=(S,\le,(P_g)_{g\in G},(Q_g)_{g\in G})\in\K$ and let
  $\pi:(S,\le)\longrightarrow(X,\le)$ be an order-isomorphism. Then we can
  find $P'_g,Q'_g\subseteq\overline{X}$ ($g\in G$) such that
  $\X=(X,\le,(P'_g)_{g\in G},(Q'_g)_{g\in G})\in\K$ and $\pi$ maps the colours
  as prescribed by $h$.
\end{itemize}
\end{lemma}

\begin{proof}
(a) Observe that $t_0\in P_h$ in $\L_h$ and $t_0\in P_{kh}$ in
$\L_{kh}$. Now continue by induction through $(T,\le)$ and construction of
$\L_h$, $\L_{kh}$.

(b) Put $P'_{gh}=P_g^\pi$ and $Q'_{gh}=Q_g^\pi$ ($g\in G$) to obtain the
result.
\end{proof}

If $(C,\le)$ is a chain and $a,b\in C$ are such that $a<b$ and there is no
$c\in C$ with $a<c<b$, we call the pair $(a,b)$ a \emph{gap} of $C$, denoted
$(a|b)$.

Let $\S=(S,\le,(P_g)_{g\in G},(Q_g)_{g\in G})\in\K$, let $h\in G$ and let 
$\L_{h}=
%%%%
\linebreak
(L_{h},\le,(P'_g)_{g\in G},(Q'_g)_{g\in G})$. An embedding
$\varphi:(L_h,\le)\longrightarrow(\overline{S},\le)$ is called \emph{nice}, if
it satisfies:
\begin{enumerate}
\item $\varphi$ preserves all suprema and infima, i.e.~$(\sup A)^\varphi=\sup
  A^\varphi$ and $(\inf A)^\varphi=\inf A^\varphi$ in $(\overline{S},\le)$ for
  each non-empty subset $A\subseteq L_h$;
\item $L_h^\varphi\subseteq\overline{S}\setminus S$ and
  $(P'_g)^\varphi\subseteq P_g$, $(Q'_g)^\varphi\subseteq Q_g$ for each $g\in
  G$;
\item whenever $x<y$ form a gap in $L_h$, then
  $\cot(x^\varphi,y^\varphi)=\aleph_0$ in $(\overline{S},\le)$;
\item $\cot(x)=\aleph_0$ if $x$ is the smallest or the largest element of
  $(L_h,\le)$.
\end{enumerate}

A structure $\S=(S,\le,(P_g)_{g\in G},(Q_g)_{g\in G})\in\K$ is called
\emph{nice}, if for each $h\in G$ and $x\in P_h$ in $\S$ there exists a nice
embedding $\varphi:\L_h\longrightarrow\S$ with $t_0^\varphi=x$. We denote
by $\Knice$ the class of all nice structures in $\K$.

\begin{proposition}\label{prop 3.4}
Let $\S=(S,\le,(P_g)_{g\in G},(Q_g)_{g\in G})\in\Knice$ be $G$-homogeneous for
$S$-intervals. Then $\S$ is $G$-homogeneous for $P$-intervals.
\end{proposition}

\begin{proof}
Choose $g,h\in G$ and $u,v\in P_g$, $x,y\in P_{gh}$ with $u<v$ and $x<y$. By
assumption, there are nice embeddings $\varphi,\varphi':\L_g\longrightarrow\S$
and $\psi,\psi':\L_{gh}\longrightarrow\S$ with $t_0^\varphi=v$,
$t_0^{\varphi'}=u$, $t_0^\psi=y$, $t_0^{\psi'}=x$.

First, from $\varphi',\psi'$ we obtain two continuously decreasing sequences
$(a_i)_{i\in\kappa}\subseteq(u,v)$ and $(b_i)_{i\in\kappa}\subseteq(x,y)$ such
that:
\begin{itemize}
\item[(i)]   $u=\inf_{i\in\kappa}a_i$ and $x=\inf_{i\in\kappa}b_i$;
\item[(ii)]  $a_i,b_i\in\overline{S}\setminus(S\cup P\cup Q)$ for each
  $i\in\kappa$;
\item[(iii)] $\cot(a_{i+1},a_i)=\aleph_0=\cot(b_{i+1},b_i)$ for each
  $i\in\kappa$.
\end{itemize}
By Lemma~\ref{lemma 3.1}, for each $i\in\kappa$ there is an isomorphism
$f_i:(a_{i+1},a_i)_S\longrightarrow(b_{i+1},b_i)_S$ permuting the colours as
prescribed by $h$. Patching these isomorphisms $f_i$ ($i\in\kappa$) together,
we obtain an isomorphism $f:[u,a_0)\longrightarrow[x,b_0)$ which permutes
the colours as prescribed by $h$.

Secondly, let $m=\min(T,\le)$. Then $m^\varphi<v$ and $m^\psi<y$, and we may
assume that $a_0<m^\varphi$ and $b_0<m^\psi$ (otherwise consider appropriate
upper segments of $\L_g^\varphi$ and $\L_{gh}^\psi$ subsequently). Observing
Lemma~\ref{lemma 3.3}, we see that $\varphi^{-1}\psi$ maps
$[m^\varphi,v]\cap L_g^\varphi$ onto $[m^\psi,y]\cap L_{gh}^\psi$ permuting the
colours in these subsets of $\S$ as prescribed by $h$. We have
$[m^\varphi,v]\setminus L_g^\varphi=\bigcup(a^\varphi,b^\varphi)$ and
$[m^\psi,y]\setminus L_{gh}^\psi=\bigcup(a^\psi,b^\psi)$ where the two unions
are taken over all gaps $(a|b)$ in $(T,\le)$. Moreover, for each gap $(a|b)$
in $(T,\le)$, we have
$\cot(a^\varphi,b^\varphi)=\aleph_0=\cot(a^\psi,b^\psi)$, and by
Lemma~\ref{lemma 3.1} there is an isomorphism
$\rho_{(a|b)}:(a^\varphi,b^\varphi)\longrightarrow(a^\psi,b^\psi)$ permuting
the colours as prescribed by $h$. Patching all these isomorphism $\rho_{(a|b)}$
together with $\varphi^{-1}\psi$ above, we obtain an isomorphism
$f':[m^\varphi,v]\longrightarrow[m^\psi,y]$ which permutes the colours as
prescribed by $h$.

Again by Lemma~\ref{lemma 3.1}, there is also such an isomorphism
$f'':(a_0,m^\varphi)\longrightarrow(b_0,m^\psi)$. Now $f\cup f''\cup f'$ maps
$[u,v]_S$ isomorphically onto $[x,y]_S$ and permutes the colours as prescribed
by $h$.
\end{proof}

\section{Construction of doubly homogeneous chains}
\label{sec 4}

In this section, we wish to prove Theorem~\ref{theorem main} first without
considering $H$, i.e.~for $H=\{1\}$. Afterwards, we will point out how to
change our constructions in order to accommodate arbitrary $\ell$-groups $H$.

Therefore, until Theorem~\ref{theorem 4.3} we make the \emph{general
  assumption} of
Section~\ref{sec 2}. We will construct a structure $\S\in\K$ which is
$G$-homogeneous for $S$-intervals. By Proposition~\ref{prop 3.4} it will
follow that $\S$ is doubly homogeneous in the $P$-colours as needed for
Theorem~\ref{theorem 2.4}.

The building blocks of our chain $(S,\le)$ are the following orderings, which
were already defined and used in \cite{DS}.

\begin{definition}[\protect{\cite{DS}}]
A chain $(L,\le)$ is called a \emph{good $\lambda$-set} if the following
conditions are satisfied:
\begin{enumerate}
\item $\lvert L\rvert=\lambda$ and $(L,\le)$ is dense and unbounded;
\item each $x\in L$ has countable coterminality;
\item whenever $x,y\in L$ with $x<y$, there is a set
  $A\subseteq[x,y]_{\overline{L}}\setminus L$ such that $\lvert
  A\rvert=\lambda$ and each $a\in A$ has countable coterminality.
\end{enumerate}
\end{definition}

We let $\K^\lambda$ comprise all structures $\S=(S,\le,(P_g)_{g\in G},
(Q_g)_{g\in G})\in\K$ for which $(S,\le)$ is a
good $\lambda$-set, and $|P_g|=|Q_g|=\lambda$ for each $g\in G$. We let
$\Knice^\lambda$ consist of all nice 
$\S\in\K^\lambda$. The following clarifies the existence of good
$\lambda$-sets.

\begin{lemma}[\protect{\cite[Lemma 4.2]{DS}}]\label{lemma 4.2}
There exists a good $\lambda$-set $(L,\le)$ of countable coterminality such
that $\chara(\overline{L})=\{(\mu,\aleph_0):\ \aleph_0\le\mu\le\lambda, \text{
  $\mu$ regular}\}$.
\end{lemma}

Here the final statement on $\chara(\overline{L})$ follows easily from the
construction for \cite[Lemma 4.2]{DS}.

The proof of Theorem~\ref{theorem main} uses two basic construction techniques
from \cite{DS} which we now describe.

\subsubsection*{Basic Construction A (defining a colour-permuting isomorphism)}

\begin{figure}[htbp]
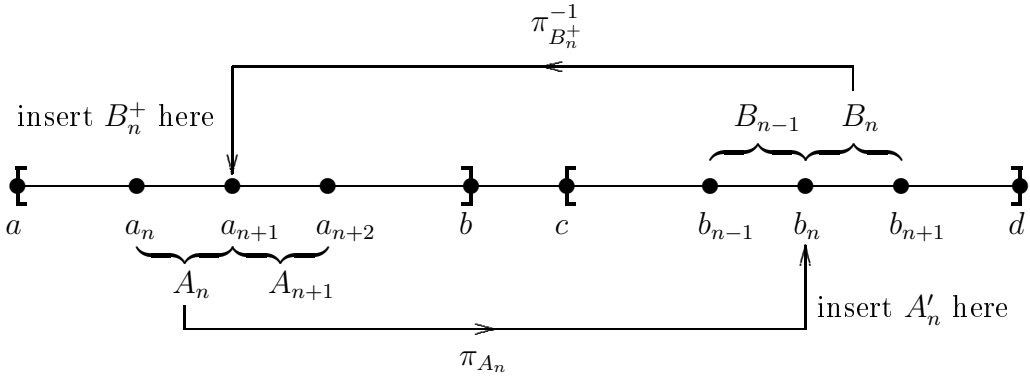

\label{fig 1}
  \begin{center}
\font\thinlinefont=cmr5
\begingroup\makeatletter\ifx\SetFigFont\undefined%
\gdef\SetFigFont#1#2#3#4#5{%
  \reset@font\fontsize{#1}{#2pt}%
  \fontfamily{#3}\fontseries{#4}\fontshape{#5}%
  \selectfont}%
\fi\endgroup%
\mbox{\beginpicture
\setcoordinatesystem units <1.00000cm,1.00000cm>
\unitlength=1.00000cm
\linethickness=1pt
\setplotsymbol ({\makebox(0,0)[l]{\tencirc\symbol{'160}}})
\setshadesymbol ({\thinlinefont .})
\setlinear
%
% Fig ELLIPSE
%
\linethickness= 0.500pt
\setplotsymbol ({\thinlinefont .})
\put{\makebox(0,0)[l]{\circle*{ 0.212}}} at 15.875 20.955
%
% Fig ELLIPSE
%
\linethickness= 0.500pt
\setplotsymbol ({\thinlinefont .})
\put{\makebox(0,0)[l]{\circle*{ 0.212}}} at 14.605 20.955
%
% Fig ELLIPSE
%
\linethickness= 0.500pt
\setplotsymbol ({\thinlinefont .})
\put{\makebox(0,0)[l]{\circle*{ 0.212}}} at 13.335 20.955
%
% Fig ELLIPSE
%
\linethickness= 0.500pt
\setplotsymbol ({\thinlinefont .})
\put{\makebox(0,0)[l]{\circle*{ 0.212}}} at 11.430 20.955
%
% Fig ELLIPSE
%
\linethickness= 0.500pt
\setplotsymbol ({\thinlinefont .})
\put{\makebox(0,0)[l]{\circle*{ 0.212}}} at 10.160 20.955
%
% Fig ELLIPSE
%
\linethickness= 0.500pt
\setplotsymbol ({\thinlinefont .})
\put{\makebox(0,0)[l]{\circle*{ 0.212}}} at  8.255 20.955
%
% Fig ELLIPSE
%
\linethickness= 0.500pt
\setplotsymbol ({\thinlinefont .})
\put{\makebox(0,0)[l]{\circle*{ 0.212}}} at  6.985 20.955
%
% Fig ELLIPSE
%
\linethickness= 0.500pt
\setplotsymbol ({\thinlinefont .})
\put{\makebox(0,0)[l]{\circle*{ 0.212}}} at  5.715 20.955
%
% Fig ELLIPSE
%
\linethickness= 0.500pt
\setplotsymbol ({\thinlinefont .})
\put{\makebox(0,0)[l]{\circle*{ 0.212}}} at 17.462 20.955
%
% Fig ELLIPSE
%
\linethickness= 0.500pt
\setplotsymbol ({\thinlinefont .})
\put{\makebox(0,0)[l]{\circle*{ 0.212}}} at  4.128 20.955
%
% Fig POLYLINE object
%
\linethickness=1pt
\setplotsymbol ({\makebox(0,0)[l]{\tencirc\symbol{'160}}})
\putrule from 11.549 21.196 to 11.412 21.196
\putrule from 11.430 21.196 to 11.430 20.698
\putrule from 11.430 20.716 to 11.549 20.716
%
% Fig POLYLINE object
%
\linethickness=1pt
\setplotsymbol ({\makebox(0,0)[l]{\tencirc\symbol{'160}}})
\putrule from 10.041 21.196 to 10.178 21.196
\putrule from 10.160 21.196 to 10.160 20.698
\putrule from 10.160 20.716 to 10.041 20.716
%
% Fig POLYLINE object
%
\linethickness=1pt
\setplotsymbol ({\makebox(0,0)[l]{\tencirc\symbol{'160}}})
\putrule from 17.344 21.196 to 17.480 21.196
\putrule from 17.462 21.196 to 17.462 20.698
\putrule from 17.462 20.716 to 17.344 20.716
%
% Fig POLYLINE object
%
\linethickness=1pt
\setplotsymbol ({\makebox(0,0)[l]{\tencirc\symbol{'160}}})
\putrule from  4.246 21.196 to  4.110 21.196
\putrule from  4.128 21.196 to  4.128 20.698
\putrule from  4.128 20.716 to  4.246 20.716
%
% Fig POLYLINE object
%
\linethickness= 0.500pt
\setplotsymbol ({\thinlinefont .})
\putrule from  4.128 20.955 to 17.462 20.955
%
% Fig POLYLINE object
%
\linethickness= 0.500pt
\setplotsymbol ({\thinlinefont .})
\putrule from 15.240 22.225 to 15.240 22.543
\putrule from 15.240 22.543 to 11.113 22.543
%
% arrow head
%
\plot 11.366 22.606 11.113 22.543 11.366 22.479 /
%
% Fig POLYLINE object
%
\linethickness= 0.500pt
\setplotsymbol ({\thinlinefont .})
\putrule from 11.113 22.543 to  6.985 22.543
\putrule from  6.985 22.543 to  6.985 21.114
%
% arrow head
%
\plot  6.921 21.368  6.985 21.114  7.048 21.368 /
%
% Fig POLYLINE object
%
\linethickness= 0.500pt
\setplotsymbol ({\thinlinefont .})
\putrule from  6.350 19.367 to  6.350 19.050
\putrule from  6.350 19.050 to 10.478 19.050
%
% arrow head
%
\plot 10.224 18.986 10.478 19.050 10.224 19.114 /
%
% Fig POLYLINE object
%
\linethickness= 0.500pt
\setplotsymbol ({\thinlinefont .})
\putrule from 10.478 19.050 to 14.605 19.050
\putrule from 14.605 19.050 to 14.605 20.161
%
% arrow head
%
\plot 14.669 19.907 14.605 20.161 14.541 19.907 /
%
% Fig TEXT object
%
\put{\SetFigFont{12}{14.4}{\familydefault}{\mddefault}{\updefault}
$a_n$
} [lB] at  5.556 20.320
%
% Fig TEXT object
%
\put{\SetFigFont{12}{14.4}{\familydefault}{\mddefault}{\updefault}
$a_{n+1}$
} [lB] at  6.826 20.320
%
% Fig TEXT object
%
\put{\SetFigFont{12}{14.4}{\familydefault}{\mddefault}{\updefault}
$a_{n+2}$
} [lB] at  8.096 20.320
%
% Fig TEXT object
%
\put{\SetFigFont{12}{14.4}{\familydefault}{\mddefault}{\updefault}
$b$
} [lB] at 10.001 20.320
%
% Fig TEXT object
%
\put{\SetFigFont{12}{14.4}{\familydefault}{\mddefault}{\updefault}
$c$
} [lB] at 11.271 20.320
%
% Fig TEXT object
%
\put{\SetFigFont{12}{14.4}{\familydefault}{\mddefault}{\updefault}
$b_{n-1}$
} [lB] at 13.176 20.320
%
% Fig TEXT object
%
\put{\SetFigFont{12}{14.4}{\familydefault}{\mddefault}{\updefault}
$b_n$
} [lB] at 14.446 20.320
%
% Fig TEXT object
%
\put{\SetFigFont{12}{14.4}{\familydefault}{\mddefault}{\updefault}
$b_{n+1}$
} [lB] at 15.716 20.320
%
% Fig TEXT object
%
\put{\SetFigFont{12}{14.4}{\familydefault}{\mddefault}{\updefault}
$\underbrace{\hspace*{0.5in}}$
} [lB] at  5.715 20.161
%
% Fig TEXT object
%
\put{\SetFigFont{12}{14.4}{\familydefault}{\mddefault}{\updefault}
$\underbrace{\hspace*{0.5in}}$
} [lB] at  6.985 20.161
%
% Fig TEXT object
%
\put{\SetFigFont{12}{14.4}{\familydefault}{\mddefault}{\updefault}
$A_n$
} [lB] at  6.191 19.526
%
% Fig TEXT object
%
\put{\SetFigFont{12}{14.4}{\familydefault}{\mddefault}{\updefault}
$A_{n+1}$
} [lB] at  7.461 19.526
%
% Fig TEXT object
%
\put{\SetFigFont{12}{14.4}{\familydefault}{\mddefault}{\updefault}
$\overbrace{\hspace*{0.5in}}$
} [lB] at 13.335 21.273
%
% Fig TEXT object
%
\put{\SetFigFont{12}{14.4}{\familydefault}{\mddefault}{\updefault}
$d$
} [lB] at 17.304 20.320
%
% Fig TEXT object
%
\put{\SetFigFont{12}{14.4}{\familydefault}{\mddefault}{\updefault}
$a$
} [lB] at  3.969 20.320
%
% Fig TEXT object
%
\put{\SetFigFont{12}{14.4}{\familydefault}{\mddefault}{\updefault}
$\overbrace{\hspace*{0.5in}}$
} [lB] at 14.605 21.273
%
% Fig TEXT object
%
\put{\SetFigFont{12}{14.4}{\familydefault}{\mddefault}{\updefault}
$B_{n-1}$
} [lB] at 13.652 21.749
%
% Fig TEXT object
%
\put{\SetFigFont{12}{14.4}{\familydefault}{\mddefault}{\updefault}
$B_n$
} [lB] at 15.081 21.749
%
% Fig TEXT object
%
\put{\SetFigFont{12}{14.4}{\familydefault}{\mddefault}{\updefault}
insert $B^+_n$ here
} [lB] at  4.128 21.749
%
% Fig TEXT object
%
\put{\SetFigFont{12}{14.4}{\familydefault}{\mddefault}{\updefault}
insert $A'_n$ here
} [lB] at 14.764 19.209
%
% Fig TEXT object
%
\put{\SetFigFont{12}{14.4}{\familydefault}{\mddefault}{\updefault}
$\pi^{-1}_{B^+_n}$
} [lB] at 10.954 23.019
%
% Fig TEXT object
%
\put{\SetFigFont{12}{14.4}{\familydefault}{\mddefault}{\updefault}
$\pi_{A_n}$
} [lB] at 10.001 18.574
\linethickness=0pt
\putrectangle corners at  3.969 23.349 and 17.585 18.451
\endpicture}
\end{center}
\caption{Defining a colour-permuting isomorphism.}
\end{figure}

Let $\S=(S,\le,(P_g)_{g\in G},(Q_g)_{g\in G})\in\K^\lambda$. Let $h\in G$ and
$a,b,c,d\in S$ with $a<b<c<d$. We will enlarge $\S$ to a superstructure
$\S^*\supseteq\S$ with $\S^*=
%%%%
\linebreak
(S^*,\le,(P_g^*)_{g\in G},(Q_g^*)_{g\in
  G})$ such that there is an isomorphism
$f:[a,b]_{S^*}\longrightarrow[c,d]_{S^*}$ permuting the colours of $\S^*$ as
prescribed by $h$. This will be obtained by splitting both $(a,b)_S$ and
$(c,d)_S$ into countably many subintervals, inserting copies of these intervals
into $(c,d)_S$ and $(a,b)_S$, respectively, to obtain $S^*$, hereby
changing the colours as prescribed by $h$, and defining the isomorphism $f$
correspondingly.

First, we choose $\mathbb{Z}$-sequences
$(a_i)_{i\in\mathbb{Z}}\subseteq[a,b]_{\overline{S}\setminus S}$ and
$(b_i)_{i\in\mathbb{Z}}\subseteq[c,d]_{\overline{S}\setminus S}$ such that
$\cot(a_i)=\cot(b_i)=\aleph_0$, $a_i<a_{i+1}$ and $b_i<b_{i+1}$ for each
$i\in\mathbb{Z}$ and $a=\inf_{i\in\mathbb{Z}}a_i$,
$b=\sup_{i\in\mathbb{Z}}a_i$, $c=\inf_{i\in\mathbb{Z}}b_i$,
$d=\sup_{i\in\mathbb{Z}}b_i$. For each $i\in\mathbb{Z}$ let $A_i'$ be a copy
of $A_i:=(a_i,a_{i+1})_S$, let $B_i^+$ be a copy of $B_i:=(b_i,b_{i+1})_S$,
and let $\pi_{A_i}:A_i\longrightarrow A_i'$, $\pi_{B_i^+}:B_i^+\longrightarrow
B_i$ be isomorphisms. Put $S^*=S\cup\bigcup_{i\in\mathbb{Z}}(A_i'\cup B_i^+)$.

We define a linear order $\le$ on $S^*$ in the natural way so that it extends
the orders of $S$ and of each $A_i'$, $B_i^+$ and the Dedekind-completions of
these sets satisfy in $(S^*,\le)$
\[
\overline{B_{i+1}}<\overline{A_i'}<\overline{B_i} \quad\text{ and }\quad
\overline{A_i}<\overline{B_i^+}<\overline{A_{i+1}} \quad\text{ for each
  $i\in\mathbb{Z}$}.
\]
We define a mapping $f:[a,b]_{S^*}\longrightarrow[c,d]_{S^*}$ as in Figure~1
                                %\ref{fig 1}
by putting $a^f=c$, $b^f=d$, $f|_{A_i}=\pi_{A_i}$ and
$f|_{B_i^+}=\pi_{B_i^+}$. Then $f$ is an order-isomorphism.

Now define colours by putting, for each $g\in G$ and $i\in\mathbb{Z}$,
\begin{align*}
P_{gh,i}' &= ((a_i,a_{i+1})\cap P_g)^{\pi_{A_i}}\subseteq\overline{A_i'},\\
Q_{gh,i}' &= ((a_i,a_{i+1})\cap Q_g)^{\pi_{A_i}}\subseteq\overline{A_i'},\\
P_{g,i}^+ &= ((b_i,b_{i+1})\cap P_{gh})^{\pi^{-1}_{B_i^+}}\subseteq\overline{B_i^+},\\
Q_{g,i}^+ &= ((b_i,b_{i+1})\cap Q_{gh})^{\pi^{-1}_{B_i^+}}\subseteq\overline{B_i^+},\\
P_g^* &= P_g\cup\bigcup_{i\in\mathbb{Z}}(P_{g,i}'\cup P_{g,i}^+),\\
Q_g^* &= Q_g\cup\bigcup_{i\in\mathbb{Z}}(Q_{g,i}'\cup Q_{g,i}^+).
\end{align*}
Then $\S^*=(S^*,\le,(P_g^*)_{g\in G},(Q_g^*)_{g\in G})\in\K^\lambda$, and it
follows that $f$ permutes the colours as prescribed by $h$. Note that if $\S$
is nice, then so is $\S^*$.\qed

\subsubsection*{Basic Construction B (extension of a colour-permuting
  isomorphism)}

Let $\U=(U,\le,(P_g)_{g\in G}, (Q_g)_{g\in G})$, $\V=(V,\le,(P'_g)_{g\in G},
(Q'_g)_{g\in G})\in \K^\lambda$ such that $\U$ is a substructure of $\V$ in the
usual sense. Assume that\\

$(*)$\label{star}\hfill
\parbox{12cm}{
whenever $v\in V\setminus S$, then the set
\[
D=\{x\in V\cup P'\cup Q':\ \forall u\in U\,:\,(u<v \longleftrightarrow u<x)\}
\]
has no greatest or smallest element and has countable coterminality, the set
$E=\{x\in U\cup P\cup Q:\ x<v\}$ has no greatest element, and the set
$F=\{x\in U\cup P\cup Q:\ v<x\}$ has no smallest element.
}\\\\

Now let $a,b,c,d\in U$ with $a<b<c<d$, let $h\in G$ and let $f\,:\,
[a,b]_U\longrightarrow [c,d]_U$ be an isomorphism permuting the colours of $\U$
as prescribed by $h$. We want to enlarge $\V$ to a superstructure $\W\supseteq
\V$ with $\W=
%%%%
\linebreak
(W,\le,(P^*_g)_{g\in G}, (Q^*_g)_{g\in G})$ such that $f$ extends
to an isomorphism $\bar{f}\,:\, [a,b]_W\longrightarrow [c,d]_W$ permuting the
colours in $\W$ as prescribed by $h$. This is achieved by inserting, for
certain decompositions $[a,b]_U=A\cup B$ with $A<B$, points into $V$ between
the sets $A$ and $B$ and also between $A^f$ and $B^f$.

Consider a decomposition $[a,b]_U=A\cup B$ with $A<B$. We distinguish between
four cases.

\paragraph{Case 1.} There is no $x\in V$ with $A<x<B$ and no $y\in V$ with
$A^f<y<B^f$.

In this case, no point is inserted, either between $A$ and $B$ or between $A^f$
and $B^f$.

\paragraph{Case 2.} There is $x\in V$ with $A< x< B$ but no $y\in V$ with $A^f<
y < B^f$.

In this case, let $X=\{x\in V:\ A< x < B\}$, let $(Y,\le)$ be a copy of
$(X,\le)$, and let $\pi\,:\, (X,\le)\longrightarrow(Y,\le)$ be an isomorphism.
We insert $Y$ into $V$ between $A^f$ and $B^f$. Next, using
Lemma~\ref{lemma 3.3}(b), we define colours in $\overline{Y}$ such that $\pi$
permutes the colours as prescribed by $h$.

\paragraph{Case 3.} There is $y\in V$ with $A^f < y < B^f$ but no $x\in V$ with
$A< x < B$.

This case is dual to Case 2.

\paragraph{Case 4.} There are $x,y\in V$ with  $A< x< B$ and $A^f < y < B^f$.

In this case, let $X=\{x\in V:\ A< x< B\}$ and $Y=\{y\in V:\ A^f < y < B^f\}$
and define $a' = \inf_{\overline{V}} X$, $b'=\sup_{\overline{V}} X$,
$c'=\inf_{\overline{V}} Y$, $d'=\sup_{\overline{V}} Y$. 
By $(*)$, $X$ and $Y$ contain no greatest or smallest element and have
countable coterminality. Hence we can deal with the intervals $[a',b']_V$ and
$[c',d']_V$ precisely as we dealt with the intervals $[a,b]_S$ and $[c,d]_S$ in
Basic Construction~A, the only difference being that the endpoints
$a',b',c',d'$ all belong to $\overline{V}$. By $(*)$, in fact we have
$a',b',c',d'\in\overline{V}\setminus(V\cup P'\cup Q')$. Hence, using
Basic Construction A we enlarge the intervals $[a',b']_V$ and
$[c',d']_V$ to intervals $[a',b']_W$ and $[c',d']_W$, respectively, define
colours in them appropriately and obtain an isomorphism
$\pi\,:\,[a',b']_W\longrightarrow [c',d']_W$ which permutes the colours as
prescribed by $h$.

If these constructions are carried out for each decomposition $[a,b]_U=A\cup
B$ with $A<B$, we clearly obtain the required extension $\W$ and, by patching
together all the ``local'' isomorphisms, also the isomorphism $\overline{f}$
as required. Note that if $\U,\V\in\Knice$, then also $\W\in\Knice$. Moreover,
we have $\V\in \K^\lambda$, and if there are at most $\lambda$ decompositions
$[a,b]_U=A\cup B$ as described, then again $\W\in\K^\lambda$. \qed

\begin{theorem}\label{theorem 4.3}
Let $G$ be any group and $\lambda\ge\lvert G\rvert$ a regular uncountable
cardinal. Then there exists $\S\in\Knice^\lambda$ with the following
properties:
\begin{enumerate}
\item $\S$ is doubly homogeneous in the $P$-colours;
\item $(S,\le)$ has countable coterminality;
\item $(\overline{S},\le)$ is not anti-isomorphic to itself.
\end{enumerate}
\end{theorem}

\begin{proof}
We will first consider the case $\aleph_0<\kappa<\lambda$ (although the case
$\aleph_0=\kappa<\lambda$ will be quite similar). We first construct
$\S\in\Knice^\lambda$ which is $G$-homogeneous for $S$-intervals.
Let $\L_e = (L_e,\le, (P^0_g)_{g\in G}, (Q^0_g)_{g\in G})$ be the coloured
linear ordering defined above.

Let $\P$ be the set of all gaps of $(L_e,\le)$. For each gap $(a|b)$ of
$(L_e,\le)$, choose a copy $L_{(a|b)}$ of the good $\lambda$-set given by
Lemma~\ref{lemma 4.2}. Define
\[
S_0=\bigcup\{L_{(a|b)}:\ (a|b)\in\P\}.
\]

Next we define a linear order on $L_e\cup S_0$ in the unique way so that it
extends the given orders of $L_e$ and each $L_{(a|b)}$ and so that
$a<L_{(a|b)}<b$ for each $(a|b)\in\P$. By construction of $L_e$, whenever
$c,d\in L_e$ with $c<d$, there is $(a|b)\in\P$ with $c\le a< b\le d$ and so
there is $s\in S_0$ with $a< s< b$. Therefore we regard $C:=L_e\setminus\{\max
L_e, \min L_e\}$ as a subset of $\overline{S_0}\setminus S_0$, and we put
\[
\S_0 = (S_0,\le, (P^0_g)_{g\in G}, (Q^0_g)_{g\in G}).
\]
Clearly, $\S_0$ satisfies conditions~\ref{def 2.2}(1)--(4). By
Lemma~\ref{lemma 4.2}, observe that
$\chara(\overline{L_{(a|b)}})=\{(\mu,\aleph_0):\ \aleph_0\le\mu\le\lambda,
\text{ $\mu$ regular}\}$ for each gap $(a|b)\in P$. Hence, if
$x\in\overline{L_{(a|b)}}$ with $\cof(x)=\lambda$, there is a continuously
increasing sequence
$(a_i)_{i<\lambda}\subseteq\overline{L_{(a|b)}}\subseteq\overline{S_0}$ such
that $x=\sup_{i<\lambda}a_i$ and $\coi(a_i)=\aleph_0$ for each
$i<\lambda$. However, note that since $\S_0$ satisfies conditions~\ref{def
  2.2}(3),(4), the above is impossible for points $x\in\overline{S_0}$ with
$x\in P^0$. Since $\chara_{\overline{S_0}}(Q^0)=\{y\in\overline{S_0}:\
\chara(y)=(\kappa,\lambda)\}$, it follows that $\S_0$ also satisfies
condition~\ref{def 2.2}(5) and hence $\S_0\in\Knice^\lambda$. Moreover, the
set $L_{(a|b)}$ with $b=\max L_e$ ($a=\min L_e$) is a final (initial) segment
of $\S_0$, respectively, and has countable coterminality, so
$\cot(S_0)=\aleph_0$.

We wish to obtain the required structure $\S=(S,\le,(P_g)_{g\in G},
(Q_g)_{g\in G})\in\Knice^\lambda$ as the union of a tower (indexed by
$\lambda$) of structures $\S_i=
%%%%
\linebreak
(S_i,\le,(P^i_g)_{g\in G}, (Q^i_g)_{g\in
 G})\in\Knice^\lambda$. For this, we employ the construction in the proof of
Theorem~2.11 (parts (I)--(III)) of \cite{DS}; here our description will be
less formal.

We let $M_i$ be the set of all quintuples $(h,a,b,c,d)$ such that $h\in G$,
$a,b,c,d\in S_i$ and $a<b<c<d$, and enumerate $M_i$ by a suitable subset
$\mu_i\subseteq \lambda$. If we deal during the construction with such a
quintuple $(h,a,b,c,d)\in M_i$ for the first time at step $i+1$, we employ
Basic Construction A to obtain $\S_{i+1}$ and an isomorphism from
$[a,b]_{S_{i+1}}$ onto $[c,d]_{s_{i+1}}$ permuting the colours as prescribed
by $h$. Later on, we employ Basic Construction B to extend isomorphisms
constructed at previous stages.

For limit ordinals $j$, we just put $\S_j=\bigcup_{i<j} \S_i$, in the
natural way. Since we perform the extension of a constructed isomorphism
$\lambda$ many times, for each $i<\lambda$ and $(h,a,b,c,d)\in M_i$ we finally
obtain an isomorphism $f\,:\,[a,b]_S\longrightarrow [c,d]_S$ permuting the
colours as prescribed by $h$. This shows that $\S$ is $G$-homogeneous for
$S$-intervals.

Also since Constructions A and B are carried out only at points
$x\in\overline{S_i}\setminus S_i$, i.e. ``inside'' $\overline{S_i}$, the set
$S_0$ remains unbounded above and below in each $S_j$ ($j<\lambda$) and in $S$.
So, $(S,\le)$ has countable coterminality.

We have to ensure that all structures $\S_i$ ($i<\lambda$) and $\S$ belong to
$\Knice^\lambda$. To see that each element $s\in S$ has countable
coterminality, observe that this is true in $S_0$, and that if $i<\lambda$,
then in the construction of $S_{i+1}$ new elements only get inserted at
points of $\overline{S_i}\setminus S_i$; so, in particular, each $s\in S_i$
has the same character in $S_{i+1}$ as in $S_i$.

Moreover, if $i<\lambda$ and $\varphi:\L_h\longrightarrow\S_i$ is a nice
embedding with $t_0^\varphi\in P^i_h$, then $\varphi:\L_h\longrightarrow \S_j$
should remain nice for each $i<j<\lambda$. This is the case, if all insertion
processes of Basic Construction A and B are only carried out at points
$x\in\overline{S_j}\setminus(S_j\cup L_h^\varphi)$. Indeed we have
$C\subseteq \overline{S_0}\setminus S_0$, and we declare all
points of $C$ ``forbidden points'' in the terminology of \cite{DS}. This means
that we are never allowed, later on, to perform insertion processes at the
cuts $x\in C$. This ensures that $C\subseteq \overline{S_i}\setminus S_i$ for
each $i<\lambda$, and $\id: L_e\longrightarrow \overline{S}$ is a nice
embedding. Furthermore, during the Basic Constructions A and B, to construct
$\S_{i+1}$ we insert copies of the whole intervals of $\S_i$ only into
particular points $x\in \overline{S_i}\setminus S_i$ with $\cot(x)=\aleph_0$
which are not ``forbidden'' in $\overline{S_i}$. Then we declare in the copy
$I'$ of $I$ all elements $x'\in\overline{I'}\setminus I'$ which correspond to
a forbidden point $x\in\overline{I}\setminus I$ is $\S_i$ as forbidden points
in $\S_{i+1}$ (cf.~\cite{DS}, pp.~256--258). This ensures that our
isomorphisms also preserve forbidden points, and if $\varphi\,:\,
\L_h\longrightarrow \S_i$ is a nice embedding, then $\varphi\,:\,
\L_h\longrightarrow \S_j$ remains nice for each $i<j<k$, and so is
$\varphi\,:\, \L_h\longrightarrow\S$.

Also, $\{x\in \overline{S_0}\setminus S_0:\ \cot(x)=\aleph_0\}$ is
$\lambda$-dense in $\overline{S_0}$. Now assume that $i<\lambda$ and
$x\in\overline{S_i}\setminus S_i$ with $\cot(x)=\aleph_0$. Suppose that in
the construction of $\S_{i+1}$ the chain $(X,\le)$ gets inserted into
$\overline{S_i}$ at $x$. Then by property $(*)$ of Basic Construction
B, we may assume that $(X,\le)$ has countable coterminality and the elements
$\inf X$, $\sup X$ of $\overline{S_{i+1}}$ become forbidden points. This
ensures that $\inf X$, $\sup X$ retain their countable character in each
$\overline{S_j}$ ($i<j<\lambda$) and in $\overline{S}$.

We have to show that the forbidden points do not prevent our constructions.
By induction, we may assume that if $i<j<\lambda$, then there are only
$<\lambda$ many points $x\in\overline{S_i}\setminus S_i$ into which elements
have been inserted in order to construct $S_j$. Hence at stage $j$, for any
$a<b$ in $S_i$ there are still $\lambda$ many cuts
$x\in[a,b]_{\overline{S_j}\setminus S_j}$ with $\cot(x)=\aleph_0$ into which
no element got inserted and which are not forbidden in $\S_j$, and these can
be used when we deal again with a quintuple in $M_i$. Moreover, then at stage
$j$ again only $<\lambda$ many new forbidden points are created, keeping the
above induction hypothesis. Since $\lambda$ is regular, also for limit
ordinals $j$ for any $a<b$ in $S_j$ we have $a,b\in S_i$ for some $i<j$, and
it follows that the set $\{x\in[a,b]_{\overline{S_j}\setminus S_j}:\
\cot(x)=\aleph_0\}$ has size $\lambda$.

Next we consider the characters of elements of $\overline{S}$. As noted
before, each inserted good $\lambda$-set $L_{(a|b)}$ satisfies
$\chara(L_{(a|b)})=\{(\mu,\aleph_0):\ \aleph_0\le\mu\le\lambda,
%%%%
\linebreak
 \text{ $\mu$
 regular}\}$. Also, for $L_e=\overline{L_e}\subseteq\overline{S}\setminus S$
we have (inside $\overline{S}$) $\chara(L_e)=\{(\mu,\aleph_0),(\aleph_0,\mu):\
\aleph_0\le\mu\le\lambda,\ \mu\not=\kappa, \text{ $\mu$
  regular}\}\cup\{(\lambda,\kappa),(\kappa,\lambda)\}$. So,
$\chara(\overline{S_0})=\chara(L_e)\cup\{(\kappa,\aleph_0)\}$. We want to
ensure that our construction neither destroys these characters nor adds new
ones. The first part is clear, since we insert sets with countable
coterminality only at points
$x\in\overline{S_i}\setminus S_i$ with $\cot(x)=\aleph_0$ and add a forbidden
point to the left and to the right of the inserted set, thereby preventing
further insertions of sets at these points.

To ensure the second part, we may proceed as in \cite[p.~231,
Cor.~2.5]{D}. That is, suppose $(j_i)_{i\in\mathbb{N}}$ is a countable sequence
of ordinals $<\lambda$ and $a_i,b_i\in S_{j_i+1}\setminus S_{j_i}$ with
$a_i<a_{i+1}<b_{i+1}<b_i$ for each $i\in\mathbb{N}$. Let
$j=\sup_{i\in\mathbb{N}}j_i$. Then in $S_j=\bigcup_{i\in\mathbb{N}}S_{j_i}$
we have $\sup\{a_i:\ i\in\mathbb{N}\}=\inf\{b_i:\
i\in\mathbb{N}\}=:x\in\overline{S_j}\setminus S_j$ and $\cot(x)=\aleph_0$ in
$\overline{S_j}$. In order to prevent sets getting inserted at later stages at
$x$, we declare $x$ (and any point arising in this way) as a forbidden point
of $S_j$. This ensures that also $\cot(x)=\aleph_0$ in
$\overline{S}$. Therefore $\chara(\overline{S})=\chara(\overline{S_0})$.

In particular, since $\aleph_0<\kappa<\lambda$, we have
$(\kappa,\aleph_0)\in\chara(\overline{S})$ but
$(\aleph_0,\kappa)\not\in\chara(\overline{S})$, so $(\overline{S},\le)$ is not
anti-isomorphic to itself.

By similar remarks as above about $\S_0$, we also obtain that each $\S_j$
($j<\lambda$) and $\S$ satisfy condition~\ref{def 2.2}(5).

Now $\S\in\Knice$ is $G$-homogeneous for $S$-intervals and hence, by
Proposition~\ref{prop 3.4}, also $G$-homogeneous for $P$-intervals. Since
$\cot(S)=\aleph_0$, a patching argument similar to the one used for
Lemma~\ref{lemma 3.1} shows that $\S$ is doubly homogeneous in $P$-colours.

It remains to consider the case $\lambda=\aleph_1$ and so $\kappa=\aleph_0$.
The above construction would also work to give $\S\in\Knice^\lambda$ with
(1), (2) and $\chara(\overline{S}) = \{(\aleph_0,\aleph_0),
(\aleph_1,\aleph_0), (\aleph_0,\aleph_1)\}$, but now this does not imply that
$(\overline{S},\le)$ is not anti-isomorphic to itself. To remedy this, we can
proceed as follows. Let $(D,\le)$ be the chain by inserting in the ordinal
$\omega_1+1$ between each element and its successor a copy of $\omega_1^*$,
the converse of $\omega_1$, and also adding a copy of $\omega_1^*$ to the right
of $\omega_1+1$. Then $D=\overline{D}$, and there is a unique element $z\in
D$ with $\chara(z)=(\aleph_1,\aleph_1)$. Now put $(D,\le)$ to the right of
$(L_e,\le)$ and proceed, with $L_e$ replaced by $L_e\cup D$, as above,
inserting copies of good $\aleph_1$-sets into each gap of the chain $(L_e\cup
D,\le)$ to obtain $\S_0$.
  
The construction above now produces a chain $(S,\le)$ with
$D\subseteq\overline{S}\setminus S$. In particular, $z\in D$ also satisfies
$\cot(z)=\aleph_1$ in $\overline{S_0}$. By construction, there is a
continuously increasing sequence $(a_i)_{i\in\omega_1}\subseteq
D\subseteq\overline{S_0}$ such that $z=\sup_{i\in\omega_1} a_i$ and
$\coi(a_i)=\aleph_1$ for each $i\in\omega_1$. There is also a continuously
decreasing sequence $(b_j)_{j\in\omega_1}\subseteq D\subseteq \overline{S_0}$
such that $z=\inf_{j\in\omega_1} b_j$ and $\cof(b_j)=\aleph_0$ for each
$j\in\omega_1$. But there is no element  $z'\in\overline{S_0}$ with
$\cot(z')=\aleph_1$ and these two asymmetric ascending respectively descending
approximation properties interchanged, so $(\overline{S_0},\le)$ is not
anti-isomorphic to itself. Since $z\in D\subseteq\overline{S}$ keeps these
properties in $(\overline{S},\le)$, but the construction produces no
$z^\prime\in \overline{S}$ with the interchanged properties, it follows that
$(\overline{S},\le)$ is not anti-isomorphic to itself.
\end{proof}

We note that we could have ensured that $(\overline{S},\le)$ is not
anti-isomorphic to itself easier, without analysing the ``interior'' character
of (elements of) $(\overline{S},\le)$, by constructing the set $(S,\le)$ with
uncountable cofinality and countable coinitiality. However, for use of
Theorem~4.3 in \cite{DGG} it will be essential that $S$ has countable
coterminality. Now we have:

\paragraph{Proof of Theorem~\ref{theorem main} in case $H=\{1\}$:} Immediate
by Theorems~\ref{theorem 4.3} and \ref{theorem 2.4}.\endproof

Next we wish to show how to change the above construction in order to embed an
arbitrary $\ell$-group into $A(S)$. First we recall:

\begin{proposition}[Holland \cite{H1}]\label{prop 4.4}
Let $H$ be any $\ell$-group. Then there exists a chain $(C,\le)$ with
$|C|=\max\{|H|,\aleph_0\}$ such that $H$ $\ell$-embeds into $A(C)$.
\end{proposition}

Next we want to embed $A(C)$ into $A(S)$, for some suitable doubly homogeneous
chain $(S,\le)$. The following tool uses ideas already contained in Holland
\cite[proof of Theorem 4]{H1}.

\begin{proposition}[\protect{\cite[Theorem 4.3]{D}}]\label{prop 4.5}
Let $(S,\le)$ be a doubly homogeneous chain. Let $C\subseteq
\overline{S}\setminus S$ such that $C$ contains all suprema and infima (taken
in $\overline{S}$) of bounded subsets of $C$ and $\cot(a,b)_S=\aleph_0$ for
each gap $(a,b)$ of $(C,\le)$. Then there exists an $\ell$-embedding of $A(C)$
into $A(S)$.
\end{proposition}

Now we obtain:

\paragraph{Proof of Theorem~\ref{theorem main} in the general case:}
By Proposition~\ref{prop 4.4}, we can choose a chain $(Y,\le)$ with
$|Y|=\max\{|H|,\aleph_0\}$ such that $H\subseteq A(Y)$ (as $\ell$-groups). Let
$Y\times\{1,2\}$ be ordered lexicographically, let $(X,\leq)$ be the
chain $\overline{Y\times\{1,2\}}$ with a copy of $\omega$ added to the right
of $Y\times\{1,2\}$. Observe that $\overline{Y\times\{1,2\}}$ is the same as
the chain $\overline{Y}$ with each element of $Y$ replaced by a $2$-element
chain. We identify $A(Y)$ with $A(Y\times\{1\})\subseteq A(X)$. Now perform
the construction of the above proof of Theorem~\ref{theorem 4.3}, with $L_e$
replaced by $L_e\cup X$, where $(X,\le)$ is placed to the right of $L_e$ if
$\lambda\ge\aleph_2$, and with $L_e\cup D$ replaced by $L_e\cup D\cup X$ with
$X$ placed to the right of $L_e\cup D$, if $\lambda=\aleph_1$. The points of
$X$ also become ``forbidden'' points in $\S_0$. Since $\cof(x),\coi(x)\le
|Y|<\lambda$ for each $x\in X$, condition~\ref{def 2.2}(5) is not disturbed
and we obtain $\S_0\in\Knice^\lambda$. Hence we can continue our construction
as before, obtaining $\S\in\Knice^\lambda$ such that (cf.~the proof of
Theorem~\ref{theorem 2.4}) the chain $(P_e,\le)$ is doubly homogeneous and
$\out(A(P_e))\cong G$. Since $X\cap P=\emptyset$ and by construction
$\cot((y,1),(y,2))=\aleph_0$ for each $y\in Y$, Proposition~\ref{prop 4.5} now
yields an $\ell$-embedding of $A(X)$, hence of $H$, into $A(P_e)$. \endproof

Next we wish to turn to the proof of Corollary~\ref{cor main}. Let $(S,\le)$
be an infinite chain and $H\subseteq A(S)$ a subgroup. Then $(H,\S)$ is called
a \emph{triply transitive ordered permutation group}, if whenever $A,B\subseteq
S$ each have three elements, there exists $h\in H$ with $A^h=B$. An element
$h\in H$ is called \emph{strictly positive}, if $h>e$ in $A(S)$, and
\emph{bounded}, if there are $a,b\in S$ with $a<b$ such that each $x\in S$
with $x^h\neq x$ satisfies $a<x<b$. For the proof of Corollary~\ref{cor main}
we need the following special case of McCleary \cite{Mc1}. It generalizes one
implication of Proposition~\ref{prop 2.1}.

\begin{proposition}[McCleary \protect{\cite[Main Theorem 4]{Mc1}}]
\label{prop 4.6}
Let $(H,S)$ be a triply transitive ordered permutation group such that $H$
contains a strictly positive bounded element and $(\overline{S},\le)$ is not
anti-isomorphic to itself. Then each automorphism $\varphi$ of $H$ is induced
by a unique automorphism $f\in A(\overline{S})$, which permutes the orbits of
$H$ in $\overline{S}$, such that
\[
h^\varphi = f^{-1}\cdot h\cdot f\quad\text{for each $h\in H$}.
\]
\end{proposition}

Using this and Theorem~\ref{theorem main} and \ref{theorem 4.3}, we can give
the

\paragraph{Proof of Corollary 1.2.} By the proof of Theorem \ref{theorem main},
in particular by Theorems \ref{theorem 4.3} and \ref{theorem 2.4}, there
exists $\S=(S,\le,(P_g)_{g\in G}, (Q_g)_{g\in G})\in\K^\lambda$ with the
properties of Theorem \ref{theorem 2.4} and such that $H\subseteq A(X)$ (as
$\ell$-groups, where $X=P_e$). In particular, $G\cong\out(A(X))$, the sets
$P_g$ ($g\in G$) are all the $A(X)$-orbits in $\ov{X}$ which are isomorphic to
$X$, and for each $g\in G$ there is $f_g\in A(\ov{X})$ permuting the
$A(X)$-orbits such that $X^{f_g} = P_g$. So, $A(X)^{f_g}=A(X)$ for each $g\in
G$.
By L\"owenheim-Skolem, we can find an
$\ell$-group $K$ of size $\lambda$ such that $H\subseteq K\subseteq A(X)$ (as
$\ell$-groups), $(K,X)$, just like $(A(X),X)$, is triply transitive and
contains a strictly positive bounded element, and $K^{f_g}=K$ for each $g\in
G$. Hence $f_g\in\aut(K)$, and by Proposition~\ref{prop 4.6}, $f_g$ permutes
the orbits of $H$ in $\overline{X}$. So, $X^{f_g}$ is an orbit of $H$, showing
that each $P_g$ ($g\in G$) is also an $H$-orbit. Now let $U$ be any $H$-orbit
in $\ov{X}$ with $U\cong X$. By condition (2.2.5) it follows that $U\subseteq
P$. Thus $U=P_g$ for some $g\in G$, showing that $H$ and $A(X)$ have the same orbits in
$\overline{X}$ which are isomorphic to $X$. Also, $\aut(H)$ and $\aut(A(X))$
induce the same permutations of these orbits. Thus
$\out(H)\cong\out(A(X))\cong G$.\endproof

Finally we just remark that in a similar way, by L\"owenheim-Skolem arguments
it follows that Theorem~\ref{theorem main} also holds for singular cardinals
$\lambda>|G|$, and hence so does Corollary~\ref{cor main}.

\end{document}